\documentclass[11pt, a4paper]{article}
\usepackage{amssymb}
\usepackage{amsthm}
\usepackage{amsmath}
\usepackage{latexsym}
\usepackage[T1]{fontenc}
\usepackage[latin1]{inputenc}

\usepackage{xcolor,verbatim,amsthm,amsmath,amssymb}
\def \bop {\noindent\textit{Proof}}
\def \eop {\hbox{}\nobreak\hfill
\vrule width 2mm height 2mm depth 0mm
\par \goodbreak \smallskip}

 \newtheorem{res}{Result}[section]
 \newtheorem{thm}[res]{Theorem}
 \newtheorem{ex}[res]{Example}

 \newtheorem{rem}[res]{Remark}

\newtheorem{prop}[res]{Proposition}
 \newtheorem{lem}[res]{Lemma}
 \newtheorem{cor}[res]{Corollary}
 \newtheorem{defi}[res]{Definition}

\numberwithin{equation}{section}

\newcommand{\integ}[2]{\displaystyle \int_{#1}^{#2}}
\newcommand{\dint}{\displaystyle\int}

\def\cE{\mathcal E}
\def\S{\mathcal S}
\def\H{\mathcal H}
\def\T{ [0,T]}
\def\hY{\hat Y}
\def\hZ{\hat Z}
\def\Lam{\Lambda}

\def\1{\mathbf 1}

\def\cE{\mathcal E}
\def\E{\mathbb E}

\def\R{\mathbb R}
\def\N{\mathbb N}
\def\bQ{\mathbb Q}
\def\Q{\mathcal Q}
\def\T{\mathbb [0,T]}
\def\F{\mathcal F}
\def\P{\mathbb P}

\def\al{\alpha}

\def\R{\mathbb R}
\def\bQ{\mathbb Q}
\def\Q{\mathcal Q}
\def\T{\mathbb T}
\def\F{\mathcal F}
\def\L{\mathbb L}
\def\P{\mathbb P}
\def\E{\mathbb E}
\def\1{\mathbf 1}

\def\L{\mathcal L}
\def\sig{\sigma}
\def\hf{\hat f}
\def\hS{\hat S}
\def\bg{\hat g}
\def\T{[0,T]}
\begin{document}
\date{}
\title{On backward stochastic differential equations driven by a family of It\^o's processes.}
\author{A. Berkaoui\\
Department of mathematics and statistics\\ Science college\\ Al-Imam Mohammed Ibn Saud Islamic University (IMSIU),\\ P.O. Box 90950,
Riyadh 11623, Saudi Arabia\\
Email : berkaoui@yahoo.fr\\ \\ \\ E. H. Essaky \\  Universit\'{e} Cadi Ayyad\\
 Facult\'{e} Poly-disciplinaire\\
Laboratoire de Mod\'{e}lisation et Combinatoire\\
D\'{e}partement de Math\'{e}matiques et Informatique\\ B.P. 4162, Safi, Maroc \\ Email: essaky@uca.ma\vspace*{0.2in}\\
}
\maketitle\maketitle

\begin{abstract}
We propose to study a new type of Backward stochastic differential equations, driven by a family of It\^o's processes. We prove existence and uniqueness of the solution, and investigate stability and comparison theorem.
\end{abstract}

\noindent \textbf{Keys Words:} Backward stochastic differential
equation; Family of It\^o's processes; Dynamic sublinear expectation operator; $m$-stability; Hedging claims under model uncertainty.
\medskip

\noindent \textbf{AMS Classification}\textit{(1991)}\textbf{: }60H10,
60H20.

\section{Introduction.}
\indent Originally motivated by questions arising in
stochastic control theory, the theory of backward stochastic differential
equations (BSDEs for short) has found important applications in fields as stochastic control, mathematical finance, Dynkin games and the second order PDE theory (see, for example, \cite{EPQ, PP2, PP1, BH, BHP} and the references therein).
\\
BSDEs have been introduced long time ago by J. B. Bismut \cite{Bi} both
as the equations for the adjoint process in the stochastic version of Pontryagin maximum principle as well as the model behind the
Black and Scholes formula for the pricing and hedging of options in mathematical finance. However the first published paper on nonlinear
BSDEs appeared only in 1990, by Pardoux and Peng \cite{PP1}. The classical BSDE consists of an equation of the form:
\begin{equation}\label{equa0}
 Y_t = \xi + \dint_t^T f(s,Y_s,Z_s)ds - \dint_t^T Z_s .dW_s,
\quad\quad 0\leq t\leq T.
\end{equation}
driven by a $d$-dimensional Brownian motion $W$, with a deterministic terminal time $T>0$, a generator $f:[0,T]\times \Omega\times\R\times\R^d\rightarrow\R$ and an $\F_T$-measurable terminal value $\xi$, where $({\F}_t)_{t\leq T}$ is the natural
filtration of $(W_t)_{t\leq T}$ augmented by the null sets. The solution of this equation, denoted by $eq(f,\xi)$, is a pair of adapted processes $(Y,Z)$ with values in $\R\times\R^d$ and $Y_T=\xi$. The existence and uniqueness result of Pardoux and Peng assumes the uniform Lipschitz assumption on the generator $f$ in $y$ and $z$. So it is supposed that there exists a positive constant $K$ such that:
$$
|f(s,\omega,y,z)-f(s,\omega,y',z')|\leq K(|y-y'|+|z-z'|).
$$

The proof of this Theorem is done in two steps. The first step consider the particular case where the generator $f$ does not depend on the variables $y$ and $z$. The process $M$ defined for $t\in[0,T]$ by
$$
M_t=\E\left(\xi+\int_0^T f(s)\,ds|\F_t\right),
$$
is a martingale, so by using the martingale representation theorem there exists an $\R^d$-valued integrable process $Z$ such that $M_t=M_0+\int_0^t Z_s.dW_s$. We define then the process $Y$ by
$$
Y_t=\E\left(\xi+\int_t^T f(s)\,ds|\F_t\right)=\xi+\int_t^T f(s)\,ds-\int_t^T\,Z_s.dW_s .
$$
The second step is based on a fixed point theorem: by introducing the Banach space $\H_{T,\beta}(\R^k)$ associated to the norm
$$
\|X\|_{T,\beta}=\left(\E\int_0^T\,e^{\beta s}|X_s|^2\,ds\right)^{1/2} ,
$$
we define the mapping $\Phi:\H_{T,\beta}(\R)\times\H_{T,\beta}(\R^d)\rightarrow \H_{T,\beta}(\R)\times\H_{T,\beta}(\R^d)$ by $\Phi(y,z)=(Y,Z)$ where $(Y,Z)$ is the solution of the BSDE with generator $f(s,y_s,z_s)$. Such solution exists from the first step. It is proved that $\Phi$ is a contraction for a specific value of $\beta$ and then admits a unique fixed point.

Several papers extended these results by taking a more general driving martingale or by assuming weak assumptions on the generator $f$ (Among others, see Antonelli \cite{A}, El Karoui and Huang \cite{KH}, Ma, Protter and Yong \cite{MPY}, Pardoux and Peng \cite{PP1, PP2, PP3}, Peng \cite{P1, P2} and Essaky and Hassani \cite{EH}). For example the Pardoux-Peng result can be extended easily to the new equation $eq(f,\xi,S)$:
\begin{eqnarray}
\label{eq2}
Y^S_t=\xi+\int_t^T\,f(s,Y^S_s,Z^S_s)\,ds-\int_t^T\,Z^S_s.dS_s,
\end{eqnarray}
driven by an It\^o process $S$ of the form $dS_t=\mu^S_t\,dt+\sig^S_t.dW_t$, where $\mu^S$ and $\sigma^S$ are respectively $\R^d$-valued and $\R^d\otimes\R^d$-valued predictable processes. By supposing that the matrix-valued process $\sig^S$ is invertible and that the process $(\sig^S)^{-1}\mu^S$ is uniformly bounded, we assure that $S$ has a unique martingale measure denoted by $\bQ^S$, equivalent to the physical probability $\P$. We remark that if $(Y^S,Z^S)$ is the solution of the equation $eq(f,\xi,S)$, then $(Y^S,Z^S\sig^S)$ is the solution of the equation $eq(f^S,\xi)$ where the generator $f^S$ is defined by $f^S(t,y,z)=f(t,y,z(\sig^S_t)^{-1})-z(\sig^S_t)^{-1}\mu^S_t$.

In this paper we consider a family $D$ of It\^o's processes instead of a single process $S$, such that each element of $D$ admits a unique equivalent martingale measure. We propose to solve the equation $eq(f,\xi,D)$:
\begin{eqnarray}
\label{eqq3}
\hY_t=\xi+\int_t^T\,f(s,\hY_s,\hZ_s) ds-\int_t^T\,\hZ_s.d\hS_s,
\end{eqnarray}

for which the solution is a triplet $(\hY,\hZ,\hS)$ satisfying the equation (\ref{eqq3}) and such that $\hS\in D$ and the process $\int_0^.\,\hZ_s.d\hS_s$ is a $\cE$-martingale with $\cE$ the dynamic sublinear expectation operator associated to the set of probability measures $\Q:=\{\bQ^S:\;S \in D\}$. We study existence and uniqueness of the solution, stability and comparison theorem.

This work is mainly motivated by pricing and hedging problems in Finance under model misspecification setting or the presence of some type of ambiguity.

\section{Main Result}

Let $(\Omega, {\F}, ({\F}_t)_{t\leq T}, P)$ be a stochastic
basis on which is defined a $d$-dimensional Brownian motion
$(W_t)_{t\leq T}$ such that $({\F}_t)_{t\leq T}$ is the natural
filtration of $(W_t)_{t\leq T}$ and ${\F}_0$ contains all
$P$-null sets of $\F$.  Note that $({\F}_t)_{t\leq T}$
satisfies the usual conditions, {\it{i.e.}} it is right continuous and
complete.
\medskip\\ Let us now introduce the following notations :
\medskip

\noindent
$\bullet$ ${\L}^{2}_T(\R)$ denotes the space of $\F_T-$measurable random variables $\xi$ satisfying $\E|\xi|^2<\infty$.\\
\noindent $\bullet$ ${\S}^{2}_T(\R)$ is the space of predictable processes $Y$ such that
$$
\|Y \|^2 = \E \sup_{s\leq T} |Y_s|^2 <\infty.
$$
\\
$\bullet$ ${\H}^{2}_T(\R^d)$ is the space of predictable processes $Z$ such that
$$
\|Z \|^2 = \E \dint_0^T |Z_s|^2 ds <\infty.
$$

\noindent$\bullet$ $\mathcal{P}$ denotes a set of predictable processes.

\noindent$\bullet$ $\S$ is the set of $\R^d$-valued It\^o's processes $S$ of the form $dS=\mu^S\,dt+\sigma^S\,dW$.

\noindent$\bullet$ $D$ is the set of process $S\in\S$ such that the process $\theta^S:=(\sig^S)^{-1}\mu^S$ belongs to $\mathcal{P}$ and bounded.\\

  We assume the following assumption {\bf (H)} :\\
\begin{enumerate}
\item [(1)] $\mathcal{P}$ is predictably convex: for all $X^1,X^2\in\mathcal{P}$ and a $\{0,1\}$-valued predictable process $h$, we have $X\in\mathcal{P}$ where $X_t=h_tX^1_t+(1-h_t)X^2_t$ for $t\in [0,T]$.
\item [(2)] For any predictable process $Z$ with values in $\R^d$, there exists some $S\in D$ such that
$$
ess\inf_{S'\in D}(\theta^{S'}_t. Z_t)=\theta^{S}_t. Z_t,
 \,\,\mbox{for all}\,\, t\in [0,T].
$$
\item [(3)] The terminal value $\xi\in \L_{T}^2(\R)$, the generator $f$ is uniformly $K-$Lipschitz with respect to $y$ and $z$ and $f(.,0,0)\in\H_{T}^2(\R)$.
\end{enumerate}
\vspace{0.1cm}

Let us now introduce the definition of our BSDE driven by a family of It\^o's processes.
\begin{defi}
A solution of the following BSDE $eq(f,\xi,D)$
\begin{eqnarray}
\label{eq3}
\hY_t=\xi+\int_t^T\,f(s,\hY_s,\hZ_s) ds-\int_t^T\,\hZ_s.d\hS_s,
\end{eqnarray}
is a triplet $(\hY,\hZ,\hS)$ satisfying equation (\ref{eq3}) such that $\hS\in D$ and the process $\int_0^.\,\hZ_s.d\hS_s$ is a $\cE$-martingale with $\cE$ the dynamic sublinear expectation operator associated to the set of probability measures $\Q:=\{\bQ^S:\;S \in D\}$.
\end{defi}

The main result of this paper is as follows.

\begin{thm}
\label{t1}
Under the assumption {\bf (H)}, the equation $eq(f,\xi,D)$ has a unique solution.

\end{thm}

Before proving Theorem \ref{t1}, we recall the definition of the $m$-stability property and state some intermediate results.
\begin{defi}
 A family $\Q$ of probability measures, all elements of which are equivalent to $\P$, is called multiplicativity stable ($m$-stable) if for all elements $\bQ^1, \bQ^2\in \Q$ with density processes $\Lam^1, \Lam^2$ and for all stopping time $\tau \leq T$, it holds that $\Lam_T := \Lam^1_{\tau}\Lam^2_{T}/\Lam^2_{\tau}$ is the density of some $\bQ\in\Q$.
\end{defi}

\begin{lem}
The set $\Q:=\{\bQ^S:\;S\in D\}$ is $m$-stable.

\end{lem}
\bop. First each $\bQ^S$ is defined via its Radon-Nikodym density $\Lam^S_T$, given by
$$
\Lam^S_T=\exp\left\{\int_0^T\,\theta^S_s.dW_s-\frac{1}{2}\int_0^T\,|\theta^S_s|^2\,ds\right\}.
$$
For $i=1,2$, let $S^i\in D$ with $\Lam^i_T$ the density of $\bQ^{S^i}$. Let a stopping time $\tau$ and define the probability measure $\bQ$ by its density $\Lam_T=\Lam^1_\tau\Lam^2_T/\Lam^2_\tau$ where $\Lam^i_\tau=\E(\Lam^i_T|\F_\tau)$. We shall prove that there exists some $S\in D$ such that $\bQ=\bQ^S$. For this we define the two processes $\mu$ and $\sig$ by $\mu_t=\1_{t<\tau}\mu^1_t+\1_{t\geq \tau}\mu^2_t$ and $\sig_t=\1_{t<\tau}\sig^1_t+\1_{t\geq\tau}\sig^2_t$ for $t\in [0,T]$ where $\mu^i$ and $\sig^i$ are associated to $S^i$ and define the process $S$ by $dS=\mu\,dt+\sig.dW$. We verify easily from assumption {\bf (H)}(1) by taking $h_t=\1_{\{t>\tau\}}$ that $\theta:=\sig^{-1}\mu\in\mathcal{P}$. So $S\in D$ and then $\bQ=\bQ^S$.

\eop

For the sake of clarity, we recall the following result which is a simplified
version of Proposition 3.1 in \cite{EPQ}.

\begin{prop}\label{prop3}
For a family of standard parameters $(f,\xi)$ and  $(f^{\alpha},\xi)$, with $\alpha$ from an arbitrary index set, let $(Y,Z)$ and $(Y^{\alpha}, Z^{\alpha})$ denote respectively the solution to the corresponding BSDEs $eq(f,\xi)$ and $eq(f^{\alpha},\xi)$. If there exists a parameter $\overline{\alpha}$ such that
$$
 f(t, Y_t, Z_t)=ess\inf_{\alpha} f^{\al}(t, Y_t, Z_t)= f^{\overline{\alpha}}(t, Y_t, Z_t)\,\  dP \times dt- a.e.,
 $$
then $Y_t = ess\inf_{\alpha} Y_t^{\alpha} =Y_t^{\overline{\alpha}}$ holds for all $t \leq T$, $P$-a.s..
\end{prop}

We suppose for the next two propositions that $\sig^S$ is the identity matrix for all $S\in D$. We define the dynamic sublinear expectation operator $\cE$, associated to the set of probability measures $\Q=\{\bQ^S:\;S\in D\}$, by $\cE_t(X)=ess\sup_{S\in D}\E^{\bQ^S}(X|\F_t)$.

\begin{prop}
\label{p1}
Let $g$ be a square integrable adapted process and let $(Y,Z)$ be the solution of the equation $eq(\bg,\xi)$ where $\bg(t,z)=ess\sup_{S\in D}(g_t-\mu^S_t z)$. Then for $t\in\T$
$$
Y_t=\cE_t\left(\xi+\int_t^T\,g_s\,ds\right).
$$

\end{prop}

\bop. For $S\in D$, let $(Y^S,Z^S)$ be the solution of the equation $eq(g,\xi,S)$. Then
$$
Y^S_t=\E^{\bQ^S}\left(\xi+\int_t^T\,g_s\,ds|\F_t\right).
$$
We remark also that $(Y^S,Z^S)$ is the solution of the equation $eq(g^S,\xi)$ where $g^S(t,z)=g_t-\mu^S_t z$ and then from Proposition \ref{prop3} we deduce that
$$
Y_t=ess\sup_{S\in D}Y^S_t=ess\sup_{S\in D}\E^{\bQ^S}\left(\xi+\int_t^T\,g_s\,ds|\F_t\right).
$$
Proposition \ref{p1} is then proved.
\eop

\begin{prop}\label{p2}
Under assumption \textbf{H}(2)-(3),
there exists a unique solution $(Y,Z)$ to the equation $eq(\hf,\xi)$ where $$\hf(t,y,z)=ess\sup_{S\in D}(f(t,y,z)-\mu^S_t z).$$
Moreover, we have
$$
Y_t=\cE_t\left(\xi+\int_t^T\,f(s,Y_s,Z_s)\,ds\right).
$$
\end{prop}

\bop. Along the proof, $C$ will denote a generic constant which may vary from line to line. The Proof is based on the Picard's approximation scheme. Let  $(Y^0,Z^0) =(0,0)$ and define  $(Y^{n+1},Z^{n+1})$ be the solution of the following BSDE:
\begin{equation}\label{nonlinear}
Y^{n+1}_t=\xi+\int_t^T ess\sup_{S\in D}\left[f(s,Y^{n}_s,Z^{n}_s)-\mu^S_s Z^{n}_s\right]\,ds-\int_t^T\,Z^{n+1}_s.dW_s.
\end{equation}
Let $n, m\in \N$ and $\beta >0$. Applying It\^{o}'s formula to
$(Y_{t}^{n+1} -Y_{t}^{m+1})^2 e^{\beta t}$ we get
\begin{equation}\label{yass}
\begin{array}{ll}
&(Y_{t}^{n+1} -Y_{t}^{m+1})^2 e^{\beta t} +
\integ{t}{T}e^{s\beta}|Z_s^{n+1} - Z_s^{m+1}|^2ds
+\integ{t}{T}\beta e^{\beta s}(Y_{s}^{n+1} -Y_{s}^{m+1})^2ds
\\ & =
2\integ{t}{T}e^{s\beta}(Y_{s}^{n+1} -Y_{s}^{m+1})
(\hf(s,Y_{s}^{n} ,Z_{s}^{n})-\hf(s,Y_{s}^{m},Z_{s}^{m}))ds
\\ &
\qquad-2\integ{t}{T}e^{s\beta}(Y_{s}^{n+1}
-Y_{s}^{m+1})(Z_{s}^{n+1}- Z_{s}^{m+1}).dW_{s},
\end{array}
\end{equation}
where $\hf(t,Y,Z)=f(t,Y,Z)-\mu^{\hS}_t Z$, for some ${\hS}\in D$.

From assumptions {\bf (H)}(3), it follows that for every $\varepsilon >0$,
$$
\begin{array}{ll}
& \quad 2(Y_{s}^{n+1} -Y_{s}^{m+1}) (\hf(s,Y_{s}^{n}
,Z_{s}^{n})-\hf(s,Y_{s}^{m},Z_{s}^{m}))\\ & \quad\qquad\leq
2\varepsilon\mid Y_{s}^{n+1} -Y_{s}^{m+1}\mid^2
+\frac{K^2}{\varepsilon}\mid Y_{s}^{n}-Y_{s}^{m}\mid^2+\frac{(K+C)^2}{\varepsilon}\mid Z_{s}^{n}-Z_{s}^{m}\mid^2,
\end{array}
$$
where we have used the fact that $\mu^{\hS}$ is bounded by a positive constant $C$. By taking $\varepsilon = \frac{\beta}{2}$, it follows then that
\begin{equation}\label{eqmea}
\begin{array}{ll}
&(Y_{t}^{n+1} -Y_{t}^{m+1})^2 e^{\beta t} +
\integ{t}{T}e^{s\beta}|Z_s^{n+1} - Z_s^{m+1}|^2ds
\\ & \leq
\frac{2K^2}{\beta}\integ{t}{T}e^{s\beta}\mid Y_{s}^{n}
-Y_{s}^{m}\mid^2 ds +\frac{2(K+C)^2}{\beta}\integ{t}{T}e^{s\beta}\mid Z_{s}^{n}
-Z_{s}^{m}\mid^2 ds
\\ & \qquad    -2\integ{t}{T}e^{s\beta}(Y_{s}^{n+1}
-Y_{s}^{m+1})(Z_{s}^{n+1}- Z_{s}^{m+1}).dW_{s}.
\end{array}
\end{equation}
Using a localization procedure, we have
\begin{equation}\label{equaa}
\begin{array}{ll}
& \E\integ{0}{T}e^{s\beta}|Z_s^{n+1} - Z_s^{m+1}|^2ds
\\ & \leq \frac{2(K+C)^2(T+1)}{\beta}\E\bigg[\displaystyle\sup_{s\leq T}e^{s\beta}\mid Y_{s}^{n}
-Y_{s}^{m}\mid^2  +  \integ{0}{T}e^{s\beta}\mid Z_{s}^{n}
-Z_{s}^{m}\mid^2 ds\bigg].
\end{array}
\end{equation}
It follows from (\ref{eqmea}) and Davis-Burkholder-Gundy
inequality that there exists a constant $c>0$, such that
$$
\begin{array}{ll}
&\E\displaystyle\sup_{t\leq T}(Y_{t}^{n+1} -Y_{t}^{m+1})^2 e^{\beta t}
\\ & \leq
\frac{2(K+C)^2(T+1)}{\beta}\E\bigg[\displaystyle\sup_{s\leq T}e^{s\beta}\mid Y_{s}^{n}
-Y_{s}^{m}\mid^2  +  \integ{0}{T}e^{s\beta}\mid Z_{s}^{n}
-Z_{s}^{m}\mid^2 ds\bigg]\\ & \qquad+c\E\bigg(\integ{0}{T}(e^{s\beta})^2\mid
Y_{s}^{n+1} -Y_{s}^{m+1}\mid^2\mid Z_{s}^{n+1}- Z_{s}^{m+1}\mid^2
ds\bigg)^{\frac12}
\\ &\leq \frac{2(K+C)^2(T+1)}{\beta}\E\bigg[\displaystyle\sup_{s\leq T}e^{s\beta}\mid Y_{s}^{n}
-Y_{s}^{m}\mid^2  +  \integ{0}{T}e^{s\beta}\mid Z_{s}^{n}
-Z_{s}^{m}\mid^2 ds\bigg]\\ & \qquad+\frac12\E\displaystyle\sup_{t\leq T}(Y_{t}^{n+1} -Y_{t}^{m+1})^2
e^{t\beta}
+\frac{c^2}{2}\E\integ{0}{T}e^{s\beta}\mid Z_{s}^{n+1}-
Z_{s}^{m+1}\mid^2 ds.
\end{array}
$$
By using inequality (\ref{equaa}), we get
$$
\begin{array}{ll}
&\E\sup_{t\leq T}(Y_{t}^{n+1} -Y_{t}^{m+1})^2 e^{\beta t}
\\ &\leq\frac{2(K+C)^2(T+1)}{\beta}(\frac{c^2}{2} +1)\E\bigg[\displaystyle\sup_{s\leq T}e^{s\beta}\mid Y_{s}^{n}
-Y_{s}^{m}\mid^2  +  \integ{0}{T}e^{s\beta}\mid Z_{s}^{n}
-Z_{s}^{m}\mid^2 ds\bigg]\\ & \qquad+\frac12\E\sup_{t\leq T}(Y_{t}^{n+1} -Y_{t}^{m+1})^2
e^{t\beta},
\end{array}
$$
and then
$$
\begin{array}{ll}
&\E\sup_{t\leq T}(Y_{t}^{n+1} -Y_{t}^{m+1})^2 e^{\beta t}
\\ &\leq \frac{4(K+C)^2(T+1)}{\beta}(\frac{c^2}{2} +1)\E\bigg[\displaystyle\sup_{s\leq T}e^{s\beta}\mid Y_{s}^{n}
-Y_{s}^{m}\mid^2  +  \integ{0}{T}e^{s\beta}\mid Z_{s}^{n}
-Z_{s}^{m}\mid^2 ds\bigg].
\end{array}
$$
Coming back to equation (\ref{eqmea}), we have
\begin{equation}\label{equa2}
\begin{array}{ll}
&\E\displaystyle\sup_{t\leq T}(Y_{t}^{n+1} -Y_{t}^{m+1})^2 e^{\beta t} +
\E\integ{0}{T}e^{s\beta}|Z_s^{n+1} - Z_s^{m+1}|^2ds
\\ & \leq
\frac{4(K+C)^2}{\beta} (c^2 +2)(T+1)\bigg(\E\displaystyle{\sup_{t\leq T}}(Y_{t}^{n} -Y_{t}^{m})^2
e^{t\beta}+\E\integ{0}{T}e^{s\beta}\mid Z_{s}^{n}
-Z_{s}^{m}\mid^2 ds\bigg).
\end{array}
\end{equation}
Taking
$\beta \geq 16(K+C)^2(c^2+2)(T+1)$ and $$\Gamma^{n+1, m+1} = \bigg(\E\sup_{t\leq T}(Y_{t}^{n+1} -Y_{t}^{m+1})^2
e^{t\beta}+
\E\integ{0}{T}e^{s\beta}|Z_s^{n+1} - Z_s^{m+1}|^2ds\bigg)^{\frac12},$$ it follows that, for every $n\geq m$
$$
\Gamma^{n+1, m+1}\leq \dfrac{1}{2} \Gamma^{n, m}\leq \bigg(\dfrac12\bigg)^m
\Gamma^{n-m+1, 1}.
$$
Using similar arguments as above and the fact that $f$ is $K-$Lipschitz and $f(.,0,0)\in\H^2_T(\R)$, it is not difficult to prove that there exists a positive constant C such that
$$
\Gamma^{n-m+1, 1}\leq  C.
$$
Therefore
$$
\Gamma^{n+1, m+1}\leq C\bigg(\dfrac12\bigg)^{m},
$$
and then
$$
\lim_{n, m\rightarrow +\infty}\E\sup_{t\leq T}(Y_{t}^{n}
-Y_{t}^m)^2 =0,\qquad \lim_{n, m\rightarrow +\infty}\E\integ{0}{T}\mid Z_{s}^{n}
-Z_{s}^{m}\mid^2 ds =0.
$$
Consequently, the sequence $ (Y^n,Z^n)$ converge to $(Y,Z)$ in ${\S}^{2}_T(\R)\times {\H}^{2}_T(\R^d)$.
Let $(\overline{Y}, \overline{Z})$ be the solution, which
exists according to the previous proposition, of the following BSDE
$$
\overline{Y}_t=\xi+\int_t^T ess\sup_{S\in D}\left[f(s,Y_s,Z_s)-\mu^S_s Z_s\right]\,ds-\int_t^T\,\overline{Z}_s.dW_s.
$$
It is not difficult to prove that there exists a constant $C>0$ such
that
$$
\begin{array}{ll}
&
\E\sup_{t\leq T}(Y_{t}^{n+1} -\overline{Y}_{t})^2
+\E\integ{0}{T}|Z_s^{n+1} - \overline{Z}_s|^2ds \\ & \quad \leq C \bigg[\E\sup_{t\leq
T}(Y_{t}^{n} -Y_t)^2+\E\integ{0}{T}|Z_s^{n} - Z_s|^2ds\bigg].
\end{array}
$$
Hence
$$
 \lim_{n\rightarrow +\infty}\Big[\E\displaystyle{\sup_{t\leq T}}(Y_{t}^{n+1}
-\overline{Y}_{t})^2 + \E\integ{0}{T}|Z_s^{n+1} -
\overline{Z}_s|^2ds\Big]=0.
$$
It follows that
$$
\E\sup_{s\leq T}\mid Y_s -\overline{Y}_s \mid^2 =
0,\,\,\mbox{and}\,\, \E\integ{0}{T}|Z_s - \overline{Z}_s|^2ds =0.
$$
Therefore $Y = \overline{Y}$ and $Z=\overline{Z}$, and then $(Y,Z)$ satisfies $eq(\xi,\hf)$. \\
Thanks to Proposition \ref{p1}, it follows from Equation (\ref{nonlinear}) that
\begin{eqnarray}
\label{**}}{Y^{n+1}_t=\cE_t\left(\xi+\int_t^T\,f(s,Y^{n}_s,Z^{n}_s)\,ds\right).
\end{eqnarray}
By taking the limit in (\ref{**}) and using the Fatou property of $\cE$ we obtain that
$$
Y_t\geq\cE_t\left(\xi+\int_t^T\,f(s,Y_s,Z_s)\,ds\right).
$$
Since there exists some $\hS\in D$ such that $\hf(t,Y_t,Z_t)=f(t,Y_t,Z_t)-\mu^{\hS}_t.Z_t$. So $(Y,Z,\hS)$ is the solution of $eq(f,\xi,S)$, we take the expectation with respect to $\bQ^S$ in both parts of this equation and obtain that
$$
Y_t=\E^{\bQ^S}\left(\xi+\int_t^T\,f(s,Y_s,Z_s)\,ds|\F_t\right)\leq\cE_t\left(\xi+\int_t^T\,f(s,Y_s,Z_s)\,ds\right).
$$
Therefore the second assertion is obtained. Proposition \ref{p2} is then proved.
\eop

Now we prove Theorem \ref{t1}.

\bop\,\textit{of Theorem} \ref{t1}. 
{\bf Existence of a solution:} Let $(Y,Z)$ be the solution of the equation $eq(\hf,\xi)$, where
$$
\hf(t,y,z)=ess\sup_{S\in D}\left(f(t,y,z)-\theta^S_t.z\right),
$$
and let $\hS\in D$ such that $\hf(t,Y,Z)=f(t,Y,Z)-\theta^{\hS}_t Z$. By applying Proposition \ref{p2} to the family $\{\theta^S:\;S\in D\}$ instead of the family $\{\mu^S:\;S\in D\}$ we get that
$$
Y_t=\cE_t\left(\xi+\int_t^T\,f(s,Y_s,Z_s)\,ds\right).
$$
But for $\hZ=Z(\sig^{\hS})^{-1}$ we have that
$$
Y_t=\xi+\int_t^T\,f(s,Y_s,Z_s)\,ds-\int_t^T\,\hZ_s.d\hS_s
$$

Then for all $t\in [0,T]$,
$$
\cE_t\left(\int_t^T\,\hZ_s.d\hS_s\right)=0.
$$
Therefore $(Y,\hZ,\hS)$ is a solution of the equation $eq(f,\xi,D)$.

{\bf Uniqueness of the solution:} Let $(\hY^1,\hZ^1,\hS^1)$ and $(\hY^2,\hZ^2,\hS^2)$ two solutions of the equation $eq(f,\xi,D)$.  By applying It\^o formula to the semi-martingale $e^{\beta s}(\hY^1_s-\hY^2_s)^2$ from $s=t$ to $s=T$ we obtain that
\begin{eqnarray}
\label{***}
e^{\beta t}(\delta\hY_t)^2+\beta\int_t^T\,e^{\beta s}(\delta\hY_s)^2\,ds+\int_t^T\,e^{\beta s}|\delta(\hZ_s\sig_s)|^2\,ds
\end{eqnarray}

$$=e^{\beta T}(\delta\hY_T)^2+2\int_t^T\,e^{\beta s}\delta\hY_s\, \delta f(s)\,ds-\int_t^T\,e^{\beta s}\delta\hY_s\, d(\delta M)_s,
$$

with $\delta\hY=\hY^1-\hY^2$, $\delta(\hZ\sig)=\hZ^1\sig^{\hS^1}-\hZ^2\sig^{\hS^2}$, $\delta f(s)=f(s,\hY^1,\hZ^1\sig^{\hS^1})-f(s,\hY^2,\hZ^2\sig^{\hS^2})$, $\delta M=M^1-M^2$ and $dM^i=\hZ^i.d\hS^i$ for $i=1,2$.
We need to show that the random variable $K:=\int_t^T\,e^{\beta s}\delta\hY_s\, d\delta M_s$ has a positive expected value under  a certain probability $\bQ\in\Q$. Let define $\mu_t=\1_{(\delta\hY_t\geq 0)}\mu^1_t+\1_{(\delta\hY_t< 0)}\mu^2_t$, $\sig_t=\1_{(\delta\hY_t\geq 0)}\sig^1_t+\1_{(\delta\hY_t< 0)}\sig^2_t$ and the process $S$ by $dS=\mu dt+\sig.dW$. From assumption {\bf (H)}(1) we have $S\in D$ and since $M^1$ and $M^2$ are $\bQ^S$-super martingales, then
$$
\eta_1:=\E^{\bQ^S}\left(\int_t^T\,e^{\beta s}(\delta\hY_s)_- d M^1_s\right)\leq 0,
$$
and
$$
\eta_2:=\E^{\bQ^S}\left(\int_t^T\,e^{\beta s}(\delta\hY_s)_+ dM^2_s\right)\leq 0,
$$
where $(\delta\hY_s)_+$ and $(\delta\hY_s)_-$ are respectively the positive and the negative parts of $\delta\hY_s$. We have also that
$$
\eta:=\int_t^T\,e^{\beta s}(\delta\hY_s)_+ d\delta M^1_s+\int_t^T\,e^{\beta s}(\delta\hY_s)_- d\delta M^2_s=\int_t^T\,e^{\beta s}(\delta\hY_s)Z_s dS_s,
$$
with $Z_s=\1_{\delta\hY_s\geq 0}\hZ^1_s+\1_{\delta\hY_s< 0}\hZ^2_s$. Therefore $\eta_3:=\E^{\bQ^S}(\eta)=0$ and $\E^{\bQ^S}(K)=-\eta_1-\eta_2+\eta_3\geq 0$.

Now by taking the expectation in (\ref{***}) with respect to $\bQ=\bQ^S$ and by following the same techniques as in Proposition 2.1 in \cite{EPQ} we obtain that $\delta\hY\equiv 0$ and $\delta(\hZ \sig)\equiv 0$.

\eop

\begin{rem}
It should be pointed out that our existence result hold true if we suppose that $det(\sigma^S(\sigma^S)^{tr})\neq 0, dP \times dt-$a.e. for all $S\in\S$ and the set $D$ is taken as follows $D=\{S\in\S:\;\theta^S:=\sigma^S(\sigma^S(\sigma^S)^{tr})^{-1}\mu^S\in\mathcal{P},\,\,\mbox{and bounded}\}$, where $\mathcal{P}$ satisfying assumption {\bf (H)}(1).
\end{rem}
An immediate consequence of Theorem \ref{t1} concerns the generalization of the martingale representation of a square integrable random variable.
\begin{cor}
\label{c1}
For any $\xi\in \L^2_T(\R)$, there exists a real number $x_0$, a driver $\hS\in D$ and a square integrable predictable $\R^d$-valued process $\hZ$ such that
$$
\xi=x_0+\int_0^T\,\hZ_s.d\hS_s,
$$
and the process $\left(\int_0^t\,\hZ_s.d\hS_s\right)_{t\in\T}$ is a $\cE$-martingale, i.e for all $t<u$ we have
$$
\cE_t\left(\int_0^u\,\hZ_s.d\hS_s\right)=\int_0^t\,\hZ_s.d\hS_s.
$$

\end{cor}

\begin{rem}
\label{r1}
The triplet $(\hY,\hZ,\hS)$ is the unique solution of the equation $eq(f,\xi, D)$ if and only if $\hY=Y$, $\hZ=Z(\sigma^{\hS})^{-1}$ and $ess\inf_{S\in D}(\theta^S.Z)=\theta^{\hS}.Z$ where the pair $(Y,Z)$ is the unique solution of the equation $eq(\hf,\xi)$ with $$\hf(t,y,z)=ess\sup_{S\in D}(f(t,y,z)-\theta^S_t z).$$

\end{rem}

Thanks to the previous remark we obtain comparison theorem of solutions as a direct consequence of Theorem 2.2 in \cite{EPQ}.

\begin{thm}
\label{t2}
Let $(\hY^i,\hZ^i,\hS^i)$ be the solution of the equation $eq(f^i,\xi^i, D)$ for $i=1,2$. We suppose that $\xi^1\geq \xi^2$ a.s and that $\delta_2 f_t:=f^1(t,\hY^2_t,\hZ^2{\hat \sig}^2_t)-f^2(t,\hY^2_t,\hZ^2{\hat \sig}^2_t)\geq 0$ a.s. Then for a.s any time $t$ we have $\hY^1_t\geq\hY^2_t$. Moreover if $\hY^1_t=\hY^2_t$ on a $\F_t$-measurable set $A$, then $\hY^1_s=\hY^2_s$ on $[t,T]\times A$.

\end{thm}

\bop. We have $f^1(t,\hY^2_t,\hZ^2{\hat \sig}^2_t)\geq f^2(t,\hY^2_t,\hZ^2{\hat \sig}^2_t)$, then $\hf^1(t,\hY^2_t,\hZ^2{\hat \sig}^2_t)\geq \hf^2(t,\hY^2_t,\hZ^2{\hat \sig}^2_t)$. Thanks to Theorem 2.2 in \cite{EPQ} and Remark \ref{r1} we obtain the result.

\eop

Another consequence of remark \ref{r1} concerns the explicit solution of linear BSDE.

\begin{cor}
\label{c2}
Let $(\alpha,\gamma)$ be a bounded $\R\times\R^d$-valued predictable process, $\varphi\in\H_{T}^2(\R)$ and $\xi\in\L^2_T(\R)$. Then the linear BSDE $eq(\xi,f,D)$ with $f(s,y,z)=\varphi_s+y\alpha_s+z.\gamma_s$,
$$
Y_t=\xi+\int_t^T \left(\varphi_s+Y_s\alpha_s+Z_s.\gamma_s\right)\,ds-\int_t^T\,Z_s.dS_s ,
$$
has a unique solution $(\hY,\hZ,\hS)$ such that $\hS$ is solution of the minimization problem
$$
ess\inf_{S\in D}(\theta^S.Z)=\theta^{\hS}.Z,
$$ for all predictable processes $Z$ and
$\hY$ is given by
$$
\hY_t=\Gamma_{t}^{-1}\E\left(\xi\,\Gamma_{T}+\int_t^T \varphi_s\Gamma_{s}\,ds|\F_t\right),
$$
where $\Gamma_{s}$ is the adjoint process defined for $s\geq 0$ by the forward linear SDE:
$$
d\Gamma_{s}=\Gamma_{s}\left(\alpha_s\,ds+(\sigma^{\hS}_s)^{-1}(\gamma_s-\mu^{\hS}_s).dW_s\right) ,
$$
and $\Gamma_{0}=1$.

\end{cor}

Next we illustrate previous results by an example of a geometric Brownian motion with volatility uncertainty.

\begin{ex} We consider the case $d=1$ and define the geometric Brownian motion $S$, solution of the equation
$$
dS_t=\mu S_t dt+ \sigma^S S_t dW_t,
$$
and define the family $D$ as the set of processes $S$ that satisfy $\sigma^S\in [\sigma^1,\sigma^2]$ where $\sigma^1$ and $\sigma^2$ are two positive real constants. So the equation $eq(f,\xi,D)$ has the unique solution $(\hY,\hZ,\hS)$ given by

$$\hY=Y, \;\hZ=\dfrac{Z}{\sigma(Z)},\; d\hS=\mu\hS dt+ \sigma(Z)\hS dW,$$
with $\sigma(Z)=\sigma^1\1_{(Z\geq 0)}+\sigma^2\1_{(Z< 0)}$ and the pair $(Y,Z)$ is the unique solution of the equation $eq(\hf,\xi)$ where $\hf(t,y,z)=f(t,y,z)-\dfrac{\mu}{\sigma(z)} z$.
\end{ex}

\section {Application to hedging claims under model uncertainty.}
We consider a financial market, which is composed of a riskless asset and $d$ risky assets. We suppose that the price of these $d+1$ assets is modelled as follows: $S^0\equiv 1$ and $S=(S^1,\ldots,S^d)$ is solution of the stochastic differential equation:
$$
dS_t=S_t\left(\mu_t dt+ \sigma_t.dW_t\right),
$$
where $\mu$ and $\sigma$ are respectively $\R^d$-valued and $\R^d\otimes\R^d$-valued predictable processes. By supposing that the matrix-valued process $\sig$ is invertible and that the process $\sig^{-1}\mu$ is uniformly bounded, we assure that $S$ has a unique martingale measure denoted by $\bQ^S$, equivalent to the physical probability $\P$, and therefore every contingent claim with payoff value $H$ at maturity time $T$ can be fully hedged, which means that there exists an $\R^d$-valued strategy $\phi$ and a price $\E^{\bQ^S}(H)$ such that $H=\E^{\bQ^S}(H)+\int_0^T\,\phi_s.dS_s$. In the Markovian case and for an European type option contract $H=f(S_T)$ we can express the strategy $\phi$ as follows: we define the function $u(t,x)=\E(f(S_{T-t})|S_0=x)$ and found out that $H=u(0,S_0)$ and $u$ is the solution of the partial differential equation:
$$
\partial_t u+\mu\partial_x u+\frac{1}{2}\sigma\sigma^*\partial^2_x u=0,
$$ and $u(T,x)=f(x)$. In a general setting we express the strategy $\phi$ in terms of the solution of a backward stochastic differential equation: $H=Y_T$ where $(Y,Z)$ is the solution of BSDE:
$$
Y_t=H-\int_t^T\,Z_s.dS_s.
$$

We may ask if the full hedge is still possible and what is the price if we suppose some uncertainty on the parameters $\mu$ and/or $\sigma$. More precisely we shall consider the situation where the vector valued parameter $\theta_t:=\sig_t^{-1}\mu_t$ varies in a random interval $[h_t,g_t]$ for $t\in[0,T]$. We denote $D=\{S\in\S:\;\theta\in[h,g]\,\,\mbox{a.s.}\}$.

\begin{thm}
\label{t3}
Let a contingent claim $H$ with $\cE(H)<\infty$. Then $H=\hY_T$ where $(\hY,\hZ,\hS)$ is the solution of the BSDE $eq(0,H,D)$:
$$
\hY_t=H-\int_t^T\,\hZ_s.d\hS_s,
$$
and the process $\hS$ is given by $d\hS=\sigma(\theta^0 dt + dW)$ with $\theta^0=(h^i\1_{(Z^i>0)}+g^i\1_{(Z^i\leq 0)})_{i=1\ldots d}$ and $(Y,Z)$ is the solution of the BSDE:
$$
Y_t=H-\int_t^T\,(h_s.Z^-_s+g_s.Z^+_s)dt -\int_t^T\,Z_s.dW,
$$
with $Z^{+}=((Z^{1})^{+},\ldots,(Z^{d})^{+})$ and $Z^{-}=((Z^{1})^{-},\ldots,(Z^{d})^{-})$, where $z^{+}$ and $z^{-}$ denote respectively the positive and the negative parts of $z$.
\end{thm}

\begin{rem}
\label{r2}
We refer to \cite{EPQ} among other references for more details and motivations on BSDE's and their applications in numerous domains.
\end{rem}

\end{document}